\documentclass[11pt]{article}

\usepackage{latexsym}
\usepackage{amssymb}

\newtheorem{thm}{Theorem}
\newtheorem{cor}{Corollary}

\newtheorem{rmr}{Remark}

\newtheorem{example}{Example}

\begin{document}
{
\begin{center}
{\Large\bf
Truncated moment problems for $J$-self-adjoint, $J$-skew-self-adjoint and $J$-unitary operators.}
\end{center}
\begin{center}
{\bf S.M. Zagorodnyuk}
\end{center}

\section{Introduction.}
During the past decade an increasing interest was devoted to the investigations of operators related to
a conjugation in a Hilbert space, see, e.g. ~\cite{cit_100_GP}, \cite{cit_200_GP2}, \cite{cit_300_Z}, \cite{cit_400_LZ} and references therein.
A conjugation $J$ in a Hilbert space $H$ is an {\it antilinear} operator on $H$ such that $J^2 x = x$, $x\in H$,
and
$$ (Jx,Jy)_H = (y,x)_H,\qquad x,y\in H. $$
The conjugation $J$ generates the following bilinear form:
\begin{equation}
\label{f1_1}
[x,y]_J := (x,Jy)_H,\qquad x,y\in H.
\end{equation}
A linear operator $A$ in $H$ is said to be $J$-symmetric ($J$-skew-symmetric) if
\begin{equation}
\label{f1_2}
[Ax,y]_J = [x,Ay]_J,\qquad x,y\in D(A),
\end{equation}
or, respectively,
\begin{equation}
\label{f1_3}
[Ax,y]_J = -[x,Ay]_J,\qquad x,y\in D(A).
\end{equation}
A linear operator $A$ in $H$ is said to be $J$-isometric if
\begin{equation}
\label{f1_3_1}
[Ax,Ay]_J = [x,y]_J,\qquad x,y\in D(A).
\end{equation}
A linear operator $A$ in $H$ is called $J$-self-adjoint ($J$-skew-self-adjoint, or $J$-unitary) if
\begin{equation}
\label{f1_4}
A = JA^*J,
\end{equation}
or
\begin{equation}
\label{f1_5}
A = -JA^*J,
\end{equation}
or
\begin{equation}
\label{f1_6}
A^{-1} = JA^*J,
\end{equation}
respectively.
In this paper we shall study the following three problems:
\begin{itemize}
\item{\textbf{Problem A.}} Given a finite set of complex numbers $\{ s_{n,m} \}_{n,m=0}^d$, $d\in \mathbb{N}$.
Find a $J$-self-adjoint operator $A$ in a Hilbert space $H$, and an element $x_0\in H$, such that
\begin{equation}
\label{f1_20}
[A^n x_0, A^m x_0]_J = s_{n,m},\qquad n,m=0,1,...,d.
\end{equation}

\item{\textbf{Problem B.}} Given a finite set of complex numbers $\{ s_{n,m} \}_{n,m=0}^d$, $d\in \mathbb{N}$.
Find a $J$-skew-self-adjoint operator $A$ in a Hilbert space $H$, and an element $x_0\in H$, such that
relation~(\ref{f1_20}) holds.

\item{\textbf{Problem C.}} Given a finite set of complex numbers $\{ s_{n,m} \}_{n,m=0}^d$, $d\in \mathbb{N}$.
Find a $J$-unitary operator $A$ in a Hilbert space $H$, and an element $x_0\in H$, such that
relation~(\ref{f1_20}) holds.
\end{itemize}

The problem A / B / C is said to be \textbf{the truncated moment problem for $J$-self-adjoint / $J$-skew-self-adjoint / $J$-unitary operators},
respectively.
These moment problems are analogs of the well-known truncated Hamburger and trigonometric moment problems, which are usually formulated in terms
of prescribed integrals of powers with respect to an unknown positive Borel measure, see, e.g.~\cite{cit_500_Akh}, \cite{cit_600_M},
\cite{cit_700_Z} and references therein. However, the operator statements for such problems,
close to our definitions of problems A-C, are also known, see~\cite[pp. 411,413]{cit_1000_GL}.

A solution of the moment problem A / B / C is said to be \textbf{canonical} if
$\mathop{\rm Lin}\nolimits\{ A^n x_0 \}_{n=0}^d = H$.
A solution of the moment problem A / B / C is said to be \textbf{almost canonical} if
$\mathop{\rm Lin}\nolimits\{ A^n x_0,\ J A^n x_0  \}_{n=0}^d = H$.
Our aim is to present conditions of the solvability for the problems A, B, C and to describe
some of their canonical solutions.
As a by-product, we shall obtain some extension results for $J$-skew-symmetric and $J$-isometric operators.
It should be said that under some conditions a description of $J$-skew-self-adjoint extensions of a $J$-skew-symmetric
operator was presented earlier, see~\cite{cit_1500_R} and references therein.
A description of $J$-unitary extensions of a $J$-isometric
operator, under some conditions which are different from our assumptions in this paper, was given in~\cite[Lemma 6]{cit_1700_K}.

\textbf{Notations.}
As usual, we denote by $\mathbb{R}, \mathbb{C}, \mathbb{N}, \mathbb{Z}, \mathbb{Z}_+$,
the sets of real numbers, complex numbers, positive integers, integers and non-negative integers,
respectively.
Let $m,n\in \mathbb{N}$.
The set of all complex matrices of size $(m\times n)$ we denote by $\mathbb{C}_{m\times n}$.
The set of all complex non-negative Hermitian matrices of size $(n\times n)$ we denote by $\mathbb{C}_{n\times n}^\geq$.
If $M\in \mathbb{C}_{m\times n}$ then $M^T$ denotes the transpose of $M$, and
$M^*$ denotes the complex conjugate of $M$. The identity matrix from $\mathbb{C}_{n\times n}$
we denote by $I_n$. By $\mathop{\rm Ker}\nolimits M$ we denote the kernel of $M$, i.e. all $x\in \mathbb{C}_{n\times 1}$ such
that $Mx = 0$.

By $\mathbb{C}^N$ we denote the finite-dimensional Hilbert space of complex column vectors
of size $N$ with the usual scalar product $(\vec x,\vec y)_{\mathbb{C}^N} = \sum_{j=0}^{N-1} x_j\overline{y_j}$,
for $\vec x,\vec y\in \mathbb{C}^N$,
$\vec x = (x_0,x_1,\ldots,x_{N-1})^T$, $\vec y = (y_0,y_1,\ldots,y_{N-1})^T$, $x_j,y_j\in \mathbb{C}$; $N\in \mathbb{N}$.

If H is a Hilbert space then $(\cdot,\cdot)_H$ and $\| \cdot \|_H$ mean
the scalar product and the norm in $H$, respectively.
Indices may be omitted in obvious cases.
For a linear operator $A$ in $H$, we denote by $D(A)$
its  domain, by $R(A)$ its range, and $A^*$ means the adjoint operator
if it exists. If $A$ is invertible then $A^{-1}$ means its
inverse. $\overline{A}$ means the closure of the operator, if the
operator is closable. If $A$ is bounded then $\| A \|$ denotes its
norm.
For a set $M\subseteq H$
we denote by $\overline{M}$ the closure of $M$ in the norm of $H$.
For an arbitrary set of elements $\{ x_n \}_{n\in I}$ in
$H$, we denote by $\mathop{\rm Lin}\nolimits\{ x_n \}_{n\in I}$
the set of all linear combinations of elements $x_n$,
and $\mathop{\rm span}\nolimits\{ x_n \}_{n\in I}
:= \overline{ \mathop{\rm Lin}\nolimits\{ x_n \}_{n\in I} }$.
Here $I$ is an arbitrary set of indices.
By $E_H$ we denote the identity operator in $H$, i.e. $E_H x = x$,
$x\in H$. In obvious cases we may omit the index $H$. If $H_1$ is a subspace of $H$, then $P_{H_1} =
P_{H_1}^{H}$ is an operator of the orthogonal projection on $H_1$
in $H$.
By $[H_1,H_2]$ we denote a set of all bounded linear operators, which map a Hilbert space $H_1$ into a Hilbert space $H_2$.

\section{A Hilbert space generated by a complex symmetric matrix.}

We shall need the following theorem.
\begin{thm}
\label{t2_1}
Let $\{ s_{n,m} \}_{n,m=0}^d$ be a finite set of complex numbers, $d\in \mathbb{N}$.
There exist a Hilbert space $H$, a conjugation $J$ in $H$, and elements $\{ x_n \}_{n=0}^d$ in $H$, such that
\begin{equation}
\label{f2_1}
[x_n, x_m]_J = s_{n,m},\qquad n,m=0,1,...,d;
\end{equation}
if and only if
the matrix $( s_{n,m} )_{n,m=0}^d$ is complex symmetric.

If the last conditions are satisfied, the Hilbert space $H$ may be chosen with $\dim H = d+1$.
\end{thm}
\textbf{Proof. }
The necessity follows directly from the property of the $J$-form: $[x,y]_J = [y,x]_J$.

\noindent
Let us check the necessity.
Suppose that the matrix $\Gamma := ( s_{n,m} )_{n,m=0}^d$ is complex symmetric.
By a corollary from Takagi's factorization~\cite[Corollary 4.4.6]{cit_2000_HJ}, there exists a matrix
$\Lambda = (a_{n,j})_{n,j=0}^d\in \mathbb{C}_{(d+1)\times (d+1)}$ such that
\begin{equation}
\label{f2_2}
\Gamma = \Lambda \Lambda^T.
\end{equation}
Set $H = \mathbb{C}^{d+1}$,
$\vec e_n = (\delta_{n,0}, \delta_{n,1}, \ldots, \delta_{n,d})^T$, $0\leq n\leq d$, and
$$ J \sum_{k=0}^d \alpha_k \vec e_k = \sum_{k=0}^d \overline{\alpha_k} \vec e_k,\qquad \alpha_k\in \mathbb{C}. $$
Elements $\{ x_n \}_{n=0}^d$ we define in the following way:
$$ x_n = \sum_{j=0}^d a_{n,j} \vec e_j,\qquad 0\leq j\leq d. $$
Then
$$ [x_n, x_m]_J = (x_n, Jx_m)_H = \sum_{j,k=0}^{d} a_{n,j} a_{m,k} \delta_{j,k} = \sum_{j=0}^{d} a_{n,j} a_{m,j} = s_{n,m}, $$
for $0\leq n,m\leq d$.
$\Box$

Let $\Gamma = ( s_{n,m} )_{n,m=0}^d$ be a complex symmetric matrix. According to the proof of the last theorem, we see  that the Gram matrix
$G = ((x_n,x_m)_H)_{n,m=0}^d$ of
the constructed sequence $\{ x_n \}_{n=0}^d$ is equal to
$$ G = \Lambda \Lambda^*. $$

If $\det \Gamma\not = 0$, then $\det \Lambda\not = 0$, $\det G\not = 0$, and  $\{ x_n \}_{n=0}^d$ form a linear basis in $\mathbb{C}^{d+1}$.
In this case we have:
\begin{equation}
\label{f2_5}
\mathop{\rm Lin}\nolimits\{ x_n \}_{n=0}^d = H.
\end{equation}
Conversely, if relation~(\ref{f2_5}) holds then $\det G\not = 0$ and therefore $\det \Gamma\not = 0$.

Notice that the matrix $\Lambda$ can be constructed explicitly, since Takagi's factorization can be computed, see~\cite[p.205]{cit_2000_HJ}.

Suppose that $\det \Gamma\not = 0$.
If $\Lambda_0\in \mathbb{C}_{d+1,l}$, $l\geq d+1$, is such that $\Lambda_0 \Lambda_0^T = \Gamma$, then
$$ \Lambda_0 \Lambda_0^T = \Lambda \Lambda^T, $$
and therefore
$$ \Lambda^{-1}\Lambda_0 \Lambda_0^T (\Lambda^{-1})^T = I_{d+1}. $$
Thus $C := \Lambda^{-1}\Lambda_0\in \mathbb{C}_{d+1,l}$ is such that $CC^T = I_{d+1}$.
Consequently,
\begin{equation}
\label{f2_7}
\Lambda_0 = \Lambda C,
\end{equation}
where $C\in \mathbb{C}_{d+1,l}$, $l\geq d+1$, is such that $CC^T = I_{d+1}$.
Conversely, any matrix of the form~(\ref{f2_7}) satisfy the condition: $\Lambda_0 \Lambda_0^T = \Gamma$.

\section{Necessary conditions for the solvability of moment problems A, B, C.}
\label{s_Necessary}

Let one of the the moment problems A,B,C be given with a set $\{ s_{n,m} \}_{n,m=0}^d$, $d\in \mathbb{N}$.
Suppose that the moment problem has a solution: a $J$-self-adjoint ($J$-skew-self-adjoint, or a $J$-unitary) operator $A$ in a Hilbert space $H$, and an element
$x_0\in H$, such that relation~(\ref{f1_20}) holds.
Set
$$ x_n := A^n x_0,\qquad 0\leq n\leq d. $$
Then
\begin{equation}
\label{f3_1}
[x_n, x_m]_J = s_{n,m},\qquad n,m=0,1,...,d.
\end{equation}
Set
$$ \Gamma := ( s_{n,m} )_{n,m=0}^d,\ H_0 := \mathop{\rm Lin}\nolimits\{ x_n \}_{n=0}^{d-1},\
\mathbf{H} := \mathop{\rm Lin}\nolimits\{ x_n, Jx_n \}_{n=0}^{d}, \mathbf{J} = J|_{\mathbf{H}}. $$
Observe that $\mathbf{J}$ is a conjugation in a Hilbert space $\mathbf{H}$.
By~(\ref{f3_1}) and the property of the $J$-form $[x,y]_J=[y,x]_J$ we conclude that the matrix $\Gamma$ is
complex symmetric.

\noindent
Consider the following operator $A_0$ in $\mathbf{H}$ with the domain $D(A) = H_0$:
$$ A_0 h = A h,\qquad h\in H_0.  $$
Observe that
\begin{equation}
\label{f3_5}
A_0 \sum_{n=0}^{d-1} \alpha_n x_n = \sum_{n=0}^{d-1} \alpha_n x_{n+1},\qquad \alpha_n\in \mathbb{C}.
\end{equation}
Of course, the operator $A_0$ is $\mathbf{J}$-symmetric (respectively $\mathbf{J}$-skew-symmetric, or $\mathbf{J}$-isometric).

\noindent
The operator $A_0$ is well-defined (as a restriction of $A$).
Consequently, the equality
\begin{equation}
\label{f3_11}
\sum_{n=0}^{d-1} \alpha_n x_n = 0,\qquad \alpha_n\in \mathbb{C},
\end{equation}
implies
\begin{equation}
\label{f3_14}
\sum_{n=0}^{d-1} \alpha_n x_{n+1} = 0.
\end{equation}

Firstly, suppose that $\mathbf{H} \not= \{ 0 \}$.
Let $\{ f_j \}_{j=0}^\rho$, $0\leq \rho\leq 2d+1$, be an orthonormal basis in $\mathbf{H}$ such that $\mathbf{J} f_j = f_j$, $0\leq j \leq\rho$.
Set
$$ a_{n,j} = (x_n, f_j),\quad \Lambda = (a_{n,j})_{0\leq n\leq d,\ 0\leq j\leq \rho}. $$
Let $\Lambda_1$ ($\Lambda_2$) be a matrix which consists of the first (respectively last) $d$ rows of the matrix $\Lambda$.
Observe that~(\ref{f3_11}) is equivalent to the following condition:
\begin{equation}
\label{f3_14_1}
0 = \left(
\sum_{n=0}^{d-1} \alpha_n x_n, f_j
\right)= \sum_{n=0}^{d-1} \alpha_n a_{n,j},\qquad 0\leq j\leq \rho,
\end{equation}
or, briefly
\begin{equation}
\label{f3_14_2}
(\alpha_0,\alpha_1,...,\alpha_{d-1}) \Lambda_1 = 0.
\end{equation}
On the other hand, condition~(\ref{f3_14}) is equivalent to the following condition:
\begin{equation}
\label{f3_14_3}
0 = \left(
\sum_{n=0}^{d-1} \alpha_n x_{n+1}, f_j
\right) = \sum_{n=0}^{d-1} \alpha_n a_{n+1,j},\qquad 0\leq j\leq \rho,
\end{equation}
or, briefly
\begin{equation}
\label{f3_14_4}
(\alpha_0,\alpha_1,...,\alpha_{d-1}) \Lambda_2 = 0.
\end{equation}
Thus, we get
\begin{equation}
\label{f3_15}
\mathop{\rm Ker}\nolimits \Lambda_1^T \subseteq \mathop{\rm Ker}\nolimits \Lambda_2^T.
\end{equation}
Moreover, we may write:
$$ s_{n,m} = [x_n,x_m]_J = (x_n, Jx_m) = \sum_{j,k=0}^{\rho} a_{n,j} a_{m,k} \delta_{j,k} $$
$$ = \sum_{j=0}^{\rho} a_{n,j} a_{m,j},\qquad 0\leq n,m\leq d. $$
Consequently, we obtain that
\begin{equation}
\label{f3_16}
\Gamma = \Lambda \Lambda^T.
\end{equation}

In the case $\mathbf{H} = \{ 0 \}$ we get $x_n = 0$, $0\leq n\leq d$, and therefore all $s_{n,m}$
are zeros. Then $\Gamma = 0$ and we may choose $\Lambda = \Gamma$, so that relations~(\ref{f3_15}) and~(\ref{f3_16}) holds.

Observe that
$$ [A_0 x_n, x_m ]_{\mathbf{J}} = [ x_{n+1}, x_m ]_{\mathbf{J}} = s_{n+1,m},\
[x_n, A_0 x_m ]_{\mathbf{J}} = s_{n,m+1}, $$
$$ [A_0 x_n, A_0 x_m ]_{\mathbf{J}} = s_{n+1,m+1},\ [x_n, x_m ]_{\mathbf{J}} = s_{n,m},\qquad 0\leq n,m\leq d-1. $$
In the case of the moment problem A / B / C, we conclude that
\begin{equation}
\label{f3_17}
s_{n+1,m} = s_{n,m+1},\qquad 0\leq n,m\leq d-1;
\end{equation}
or
\begin{equation}
\label{f3_18}
s_{n+1,m} = -s_{n,m+1},\qquad 0\leq n,m\leq d-1;
\end{equation}
or
\begin{equation}
\label{f3_19}
s_{n+1,m+1} = s_{n,m},\qquad 0\leq n,m\leq d-1,
\end{equation}
respectively.

From the preceding considerations we obtain the following result.
\begin{thm}
\label{t3_1}
Let one of the moment problems A,B,C be given with a set $\{ s_{n,m} \}_{n,m=0}^d$, $d\in \mathbb{N}$.
If the moment problem A (B or C) has a solution, then
the matrix $\Gamma := ( s_{n,m} )_{n,m=0}^d$ is complex symmetric,
relation~(\ref{f3_17}) holds (respectively relation (\ref{f3_18}) holds, or relation (\ref{f3_19}) holds),
and
there exists a representation of the matrix
$\Gamma$ of the form $\Gamma = \Lambda \Lambda^T$,
where $\Lambda\in \mathbb{C}_{d+1,\rho + 1}$, $0\leq\rho\leq 2d+1$, such that
$\mathop{\rm Ker}\nolimits \Lambda_1^T \subseteq \mathop{\rm Ker}\nolimits \Lambda_2^T$.
Here $\Lambda_1$ ($\Lambda_2$) is a matrix which consists of the first (respectively last) $d$ rows of the matrix $\Lambda$.
\end{thm}

\section{Sufficient conditions for the solvability of moment problems A, B, C.}
Let one of the moment problems A, B, C be given with a set $\{ s_{n,m} \}_{n,m=0}^d$, $d\in \mathbb{N}$.
Suppose that
the matrix $\Gamma := ( s_{n,m} )_{n,m=0}^d$ is complex symmetric,
relation~(\ref{f3_17}) holds (respectively relation (\ref{f3_18}) holds, or relation (\ref{f3_19}) holds),
and
there exists a representation of the matrix
$\Gamma$ of the form
\begin{equation}
\label{f4_1}
\Gamma = \Lambda \Lambda^T,
\end{equation}
where $\Lambda\in \mathbb{C}_{d+1,\rho + 1}$, $0\leq\rho\leq 2d+1$, such that
$\mathop{\rm Ker}\nolimits \Lambda_1^T \subseteq \mathop{\rm Ker}\nolimits \Lambda_2^T$.
Here $\Lambda_1$ ($\Lambda_2$) is a matrix which consists of the first (respectively last) $d$ rows of the matrix $\Lambda$.
Let $\Lambda = (a_{n,j})_{0\leq n\leq d,\ 0\leq j\leq \rho}$, $a_{n,j}\in \mathbb{C}$.

\noindent
Set $H = \mathbb{C}^{\rho + 1}$,
$\vec e_n = (\delta_{n,0}, \delta_{n,1}, \ldots, \delta_{n,\rho})^T$, $0\leq n\leq \rho$, and
$$ J \sum_{k=0}^\rho \alpha_k \vec e_k = \sum_{k=0}^\rho \overline{\alpha_k} \vec e_k,\qquad \alpha_k\in \mathbb{C}. $$
Then $J$ is a conjugation in $H$.
Elements $\{ x_n \}_{n=0}^d$ we define in the following way:
$$ x_n = \sum_{j=0}^\rho a_{n,j} \vec e_j,\qquad 0\leq n\leq d. $$
Then
$$ [x_n, x_m]_J = (x_n, Jx_m)_H = \sum_{j,k=0}^{\rho} a_{n,j} a_{m,k} \delta_{j,k} = \sum_{j=0}^{\rho} a_{n,j} a_{m,j} = s_{n,m}, $$
\begin{equation}
\label{f4_4}
0\leq n,m\leq d.
\end{equation}
Define the following operator:
\begin{equation}
\label{f4_5}
A_0 \sum_{n=0}^{d-1} \alpha_n x_n = \sum_{n=0}^{d-1} \alpha_n x_{n+1},\qquad \alpha_n\in \mathbb{C},
\end{equation}
with the domain $D(A_0) = \mathop{\rm Lin}\nolimits\{ x_n \}_{n=0}^{d-1}$.
Let us check that it is well-defined. The latter fact means that relation~(\ref{f3_11})
implies~(\ref{f3_14}). Suppose that~(\ref{f3_11}) holds.
Then
relations~(\ref{f3_14_1}) and~(\ref{f3_14_2}) hold, with $f_j = \vec e_j$.
By our assumptions we conclude that relations~(\ref{f3_14_4}) and~(\ref{f3_14_3}) hold, with $f_j = \vec e_j$.
Then relation~(\ref{f3_14}) holds. Thus, $A_0$ is well-defined.

Choose arbitrary $h = \sum_{k=0}^{d-1} \alpha_k x_k$, $g = \sum_{j=0}^{d-1} \beta_j x_j$, $\alpha_k,\beta_j\in \mathbb{C}$, from $D(A_0)$. Then
$$ [A_0 h, g]_J = \sum_{k,j=0}^{d-1} \alpha_k \beta_j [x_{k+1}, x_j]_J = \sum_{k,j=0}^{d-1} \alpha_k \beta_j s_{k+1,j}, $$
$$ [h, A_0 g]_J = \sum_{k,j=0}^{d-1} \alpha_k \beta_j [x_{k}, x_{j+1}]_J = \sum_{k,j=0}^{d-1} \alpha_k \beta_j s_{k,j+1}. $$
In the case of the moment problem A (B) by relation~(\ref{f3_17}) (respectively by relation~(\ref{f3_18})) we conclude that
$A_0$ is $J$-symmetric (respectively $J$-skew-symmetric).
In the case of the moment problem C, by a similar argument we obtain that
$A_0$ is $J$-isometric.

\begin{thm}
\label{t4_1}
Let one of the moment problems A, B, C be given with a set $\{ s_{n,m} \}_{n,m=0}^d$, $d\in \mathbb{N}$.
Suppose that
the matrix $\Gamma := ( s_{n,m} )_{n,m=0}^d$ is complex symmetric,
relation~(\ref{f3_17}) holds (respectively relation (\ref{f3_18}) holds, or relation (\ref{f3_19}) holds),
and
there exists a representation of the matrix
$\Gamma$ of the form $\Gamma = \Lambda \Lambda^T$,
where $\Lambda\in \mathbb{C}_{d+1,\rho + 1}$, $0\leq\rho\leq 2d+1$, such that
$\mathop{\rm Ker}\nolimits \Lambda_1^T \subseteq \mathop{\rm Ker}\nolimits \Lambda_2^T$.
Here $\Lambda_1$ ($\Lambda_2$) is a matrix which consists of the first (respectively last) $d$ rows of the matrix $\Lambda$.
Suppose that for the corresponding $J$-symmetric (respectively $J$-skew-symmetric, or $J$-isometric) operator $A_0$ from~(\ref{f4_5}) there exists
a $J$-self-adjoint (respectively $J$-skew-self-adjoint, or $J$-unitary) extension in a possibly larger Hilbert space (with an extension of $J$).
Then the moment problem has a solution.
\end{thm}
\textbf{Proof.} Let $A\supseteq A_0$ be a $J$-self-adjoint ($J$-skew-self-adjoint, or $J$-unitary)
extension of the $J$-symmetric (respectively $J$-skew-symmetric, or $J$-isometric) operator $A_0$ from~(\ref{f4_5}).
By the induction argument we conclude that
\begin{equation}
\label{f4_6}
A^n x_0 = x_n,\qquad 0\leq n\leq d.
\end{equation}
By~(\ref{f4_4}) and~(\ref{f4_6}) we obtain that relation~(\ref{f1_20}) holds.
$\Box$

Thus, similar to the case of classical moment problems, we arrive to a problem of an extension of the corresponding operator.

In the case $\det\Gamma \not=0$, the above sufficient conditions of the solvability may be simplified.

\begin{cor}
\label{c4_1}
Let one of the moment problems A, B, C be given with a set $\{ s_{n,m} \}_{n,m=0}^d$, $d\in \mathbb{N}$.
Suppose that
the matrix $\Gamma := ( s_{n,m} )_{n,m=0}^d$ is complex symmetric,
relation~(\ref{f3_17}) holds (respectively relation (\ref{f3_18}) holds, or relation (\ref{f3_19}) holds),
and
$\det\Gamma \not= 0$.
Suppose that for the corresponding $J$-symmetric (respectively $J$-skew-symmetric, or $J$-isometric) operator $A_0$ from~(\ref{f4_5})
(acting in $\mathbb{C}^\rho$, $\rho = d+1$, with $\Lambda$ in~(\ref{f4_1}), provided by the corollary from Takagi's factorization)
there exists
a $J$-self-adjoint (respectively $J$-skew-self-adjoint, or $J$-unitary) extension in a possibly larger Hilbert space (with an extension of $J$).
Then the moment problem has a solution.
\end{cor}
\textbf{Proof.}
By the above-mentioned corollary from Takagi's factorization~(\cite[Corollary 4.4.6]{cit_2000_HJ}), there exists a matrix
$\Lambda\in \mathbb{C}_{(d+1)\times (d+1)}$ such that
$$ \Gamma = \Lambda \Lambda^T. $$
It is clear that $\det\Lambda\not= 0$.
Then
$\mathop{\rm Ker}\nolimits \Lambda_1^T \subseteq \mathop{\rm Ker}\nolimits \Lambda_2^T$, where
$\Lambda_1$ ($\Lambda_2$) is a matrix which consists of the first (respectively last) $d$ rows of the matrix $\Lambda$.
In fact, if $(\alpha_0,\alpha_1,...,\alpha_{d-1})\Lambda_1 = 0$, $\alpha_j\in \mathbb{C}$, then
$(\alpha_0,\alpha_1,...,\alpha_{d-1},0)\Lambda = 0$. Therefore all $\alpha_j$ are zeros.

\noindent
It remains to apply Theorem~\ref{t4_1} to complete the proof.
$\Box$

Let $H_1,H_2$ be some Hilbert spaces. An operator $J$, which maps $H_1$ into $H_2$ is said to be anti-isometric, if (\cite{cit_15000_Rayh})
$$ (Jx,Jy)_{H_2} = (y,x)_{H_1},\qquad x,y\in H_1. $$
An operator $A\in [H_1,H_2]$ is said to be $J$-self-adjoint ($J$-skew-self-adjoint) if
$B = JB^* J$  (respectively $B = -JB^* J$).

\begin{thm}
\label{t4_3} (\cite[Theorem 1]{cit_15000_Rayh})
Let $B$ be a bounded closed $J$-symmetric operator in a Hilbert space $H$, with the domain $H_1 := D(B)$ and
$H_2 := H\ominus D(B) \not= \{ 0 \}$. Let $B^*$ be the adjoint to $B$, viewed as an operator from $[H_1, H]$.
The following formula:
\begin{equation}
\label{f4_19}
\widehat B = B P^H_{H_1} + (JB^*J + S) P^H_{H_2}
\end{equation}
establishes a one-to-one correspondence between a set of all bounded $J$-self-adjoint extensions on the whole $H$ of $B$, and
a set of all $J$-self-adjoint operators $S\in [H_2, JH_2]$.
\end{thm}

The following analog of Theorem~\ref{t4_3} holds.

\begin{thm}
\label{t4_3_1}
Let $B$ be a bounded closed $J$-skew-symmetric operator in a Hilbert space $H$, with the domain $H_1 := D(B)$ and
$H_2 := H\ominus D(B) \not= \{ 0 \}$. Let $B^*$ be the adjoint to $B$, viewed as an operator from $[H_1, H]$.
The following formula:
\begin{equation}
\label{f4_20}
\widehat B = B P^H_{H_1} + (-JB^*J + S) P^H_{H_2}
\end{equation}
establishes a one-to-one correspondence between a set of all bounded $J$-skew-self-adjoint extensions on the whole $H$ of $B$, and
a set of all $J$-skew-self-adjoint operators $S\in [H_2, JH_2]$.
\end{thm}
\textbf{Proof.} The proof is similar to the proof of the Theorem~\ref{t4_3}, given in~\cite{cit_15000_Rayh}.
Since $B$ is $J$-skew symmetric, it follows that
\begin{equation}
\label{f4_21}
P^H_{JH_1} B \subseteq -J B^* J,
\end{equation}
where $B^*$ is understood for $B$, which is viewed as an operator from $[H_1, H]$.
The following property holds:
\begin{equation}
\label{f4_23}
P^H_{F} J = J P^H_{JF},
\end{equation}
where $F$ is an arbitrary subspace of $H$.
Observe that
\begin{equation}
\label{f4_25}
\widehat B^* = B^* - P^H_{H_2} JB P^H_{H_1} J + S^* P^H_{JH_2}.
\end{equation}
By~(\ref{f4_25}),(\ref{f4_21}) and~(\ref{f4_23}) we obtain that $J \widehat B^* J = - \widehat B$.

On the other hand, let $\widetilde B$ be an arbitrary bounded $J$-skew-self-adjoint extension of $B$ on the whole $H$.
Let $\widehat B_0$ be the operator $\widehat B$ from~(\ref{f4_20}) with $S=0$. Since $(\widetilde B - \widehat B_0) H_1 = \{ 0 \}$, then
$R\left( (\widetilde B - \widehat B_0)^* \right)\subseteq H_2$. We have
$R\left( \widetilde B - \widehat B_0 \right) = R\left( J(\widetilde B - \widehat B_0)^*J \right)\subseteq JH_2$.
We set $\mathcal{S} = \left( \widetilde B - \widehat B_0 \right)|_{H_2}\in [H_2,JH_2]$.
Observe that $\mathcal{S}$ is $J$-skew-self-adjoint, and $\widetilde B = \widehat B_0 + SP^H_{H_2}$ has the required form.
$\Box$

Notice that $S=0$ is $J$-self-adjoint and $J$-skew-self-adjoint. Therefore, by Theorems~\ref{t3_1}, \ref{t4_1}, \ref{t4_3}, \ref{t4_3_1} and
Corollary~\ref{c4_1} we obtain
the following two theorems.

\begin{thm}
\label{t4_4}
Let the moment problem A (B) be given with a set $\{ s_{n,m} \}_{n,m=0}^d$, $d\in \mathbb{N}$.
The moment problem A (B) has a solution if and only if
the matrix $\Gamma := ( s_{n,m} )_{n,m=0}^d$ is complex symmetric,
relation~(\ref{f3_17}) (respectively~(\ref{f3_18})) holds,
and
there exists a representation of the matrix
$\Gamma$ of the form $\Gamma = \Lambda \Lambda^T$,
where $\Lambda\in \mathbb{C}_{d+1,\rho + 1}$, $0\leq\rho\leq 2d+1$, such that
$\mathop{\rm Ker}\nolimits \Lambda_1^T \subseteq \mathop{\rm Ker}\nolimits \Lambda_2^T$.
Here $\Lambda_1$ ($\Lambda_2$) is a matrix which consists of the first (respectively last) $d$ rows of the matrix $\Lambda$.
\end{thm}

\begin{thm}
\label{t4_5}
Let the moment problem A (B) be given with a set $\{ s_{n,m} \}_{n,m=0}^d$, $d\in \mathbb{N}$.
Suppose that the matrix $\Gamma := ( s_{n,m} )_{n,m=0}^d$ is complex symmetric,
relation~(\ref{f3_17}) (respectively~(\ref{f3_18})) holds,
and $\det\Gamma \not= 0$.
Then the moment problem has a solution.
\end{thm}

Moreover, if conditions of Theorem~\ref{t4_5} are satisfied, then canonical solutions in a Hilbert space $H$, constructed as in
Corollary~\ref{c4_1}, can be obtained
by formula~(\ref{f4_19}).

On the other hand, suppose that we only know that conditions of Theorem~\ref{t4_4} are satisfied. In this case,
we may proceed as at the beginning of Section~\ref{s_Necessary} and construct a $\mathbf{J}$-symmetric ($\mathbf{J}$-skew-symmetric) operator
$A_0$ in a finite-dimensional Hilbert space $\mathbf{H}$. By Theorem~\ref{t4_3} (Theorem~\ref{t4_3_1}) extending (if necessary) this operator to
a $\mathbf{J}$-self-adjoint
(respectively $\mathbf{J}$-skew-self-adjoint) operator in $\mathbf{H}$, we obtain an almost canonical solution of
the moment problem A (respectively B).
Thus, \textit{if the moment problem A (B) is solvable, then it always has an almost canonical solution}.

\section{Extensions of $J$-isometric operators. Applications to the moment problem C.}

Let $V$ be a bounded closed $J$-isometric operator in a Hilbert space $H$, with the domain $H_0 := D(V)$ and
$H_1 := H\ominus D(V) \not= \{ 0 \}$.
Let $W \supseteq V$ be a bounded operator, defined on the whole $H$. Choose arbitrary $h,g\in H$,
$h = h_0 + h_1$, $g = g_0 + g_1$, $g_j,h_j\in H_j$, $j=1,2$.
We may write:
$$ [Wh, Wg]_J = [Wh_0 + Wh_1, Wg_0 + Wg_1]_J $$
$$ = [h_0,g_0]_J + [Vh_0, Wg_1]_J + [Wh_1, Vg_0]_J + [Wh_1, Wg_1]_J; $$
$$ [h, g]_J = [h_0 + h_1, g_0 + g_1]_J = [h_0,g_0]_J + [h_0, g_1]_J + [h_1, g_0]_J + [h_1, g_1]_J. $$
Thus, $W$ is $J$-isometric if and only if
$$ [Vh_0, Wg_1]_J + [Wh_1, Vg_0]_J + [Wh_1, Wg_1]_J = [h_0, g_1]_J + [h_1, g_0]_J + [h_1, g_1]_J. $$
In the case $h_1 = 0$, we get
\begin{equation}
\label{f5_3}
[Vh_0, Wg_1]_J = [h_0, g_1]_J,\qquad h_0\in H_0,\ g_1\in H_1.
\end{equation}
In the case $g_1 = 0$, we get
\begin{equation}
\label{f5_4}
[Wh_1, Vg_0]_J = [h_1, g_0]_J,\qquad h_1\in H_1,\ g_0\in H_0.
\end{equation}
Therefore
\begin{equation}
\label{f5_5}
[Wh_1, Wg_1]_J = [h_1, g_1]_J,,\qquad h_1,g_1\in H_1.
\end{equation}
It is clear that conditions~(\ref{f5_3}) and~(\ref{f5_4}) are equivalent.
Consequently, \textit{the operator $W$ is $J$-isometric iff relation~(\ref{f5_3}) holds and the operator $W|_{H_1}$
is $J$-isometric}.

Assume that \textit{$V$ has a bounded inverse}. Denote by $(V^{-1})^*$ the adjoint to $V^{-1}$, viewed as an operator from $[R(V), H_0]$.
Rewrite condition~(\ref{f5_3}) in the following form:
$$ (Vh_0, JWg_1) = (V^{-1} V h_0, P^H_{H_0} Jg_1) = (Vh_0, (V^{-1})^*P^H_{H_0} Jg_1), $$
where $h_0\in H_0$, $g_1\in H_1$. Then
$$ (Vh_0, JWg_1 - (V^{-1})^*P^H_{H_0} Jg_1) = 0,\qquad h_0\in H_0,\ g_1\in H_1; $$
and therefore $Wg_1 - J(V^{-1})^*P^H_{H_0} Jg_1 \in H\ominus JR(V)$, $\forall g_1\in H_1$.
Then
\begin{equation}
\label{f5_7}
P^H_{JR(V)} W|_{H_1} = J(V^{-1})^*P^H_{H_0} J|_{H_1}.
\end{equation}
From our considerations we obtain the following theorem.

\begin{thm}
\label{t5_1}
Let $V$ be a bounded closed $J$-isometric operator in a Hilbert space $H$, with the domain $H_0 := D(V)$ and
$H_1 := H\ominus D(V) \not= \{ 0 \}$.
Let $V$ have a bounded inverse. The operator $V$ can be extended to
a bounded $J$-isometric operator, defined on the whole $H$, if and only if
there exists a bounded $J$-isometric operator $W_1$ in $H$, with the domain $H_1$, such that
\begin{equation}
\label{f5_9}
P^H_{JR(V)} W_1 = J(V^{-1})^*P^H_{H_0} J|_{H_1}.
\end{equation}
Here by $(V^{-1})^*$ the adjoint to $V^{-1}$, viewed as an operator from $[R(V), H_0]$, is understood.
If such an operator $W_1$ exists then the following formula:
\begin{equation}
\label{f5_10}
W (h_0 + h_1) = V h_0 + W_1 h_1,\qquad h_0\in H_0,\ h_1\in H_1,
\end{equation}
establishes a one-to-one correspondence between a set of bounded $J$-isometric extensions $W\supseteq V$, $D(W)=H$,
and a set of all
bounded $J$-isometric operators $W_1$ in $H$, with the domain $H_1$, such that
relation~(\ref{f5_9}) holds.
\end{thm}

In general, it is not easy to construct an operator $W_1$ with the properties, described in Theorem~\ref{t5_1}.
However, in the case of a finite-dimensional $H$ and $\dim H_1 = 1$, a full explicit answer on a question of the possibility
of the extension can be given.

\begin{thm}
\label{t5_2}
Let $V$ be an invertible $J$-isometric operator in a finite-dimensional Hilbert space $H$, with the domain $H_0 := D(V)$ and
$H_1 := H\ominus D(V) \not= \{ 0 \}$, $\dim H_1 = 1$.
Let $u$ be a non-zero element in $H_1$, and $v$ be a non-zero element in $H\ominus JR(V)$.
The operator $V$ can be extended to
a $J$-isometric operator, defined on the whole $H$, if and only if
there exists a solution of the following quadratic equation with respect to an unknown $\lambda\in \mathbb{C}$:
$$ [v,v]_J \lambda^2 + 2 \left( v, (V^{-1})^*P^H_{H_0} J u \right)_H \lambda $$
\begin{equation}
\label{f5_11}
+
\overline{
\left[ (V^{-1})^*P^H_{H_0} J u, (V^{-1})^*P^H_{H_0} J u \right]_J
}
- [u,u]_J = 0.
\end{equation}
Here by $(V^{-1})^*$ the adjoint to $V^{-1}$, viewed as an operator from $[R(V), H_0]$, is understood.
If equation~(\ref{f5_11}) is solvable then the following formula:
\begin{equation}
\label{f5_12}
W (h_0 + \beta u) = V h_0 + \beta \lambda v + \beta J(V^{-1})^*P^H_{H_0} J u,\qquad h_0\in H_0,\ \beta\in \mathbb{C},
\end{equation}
establishes a one-to-one correspondence between a set of solutions $\lambda$ of equation~(\ref{f5_11})
and a set of all
$J$-isometric operators $W\supseteq V$, $D(W)=H$.
All such operators $W$ are $J$-unitary operators.
\end{thm}

\begin{rmr}
\label{r5_1}
It is clear that equation~(\ref{f5_11}) has no solutions if and only if
\begin{equation}
\label{f5_14}
[v,v]_J = 0,\quad \left( v, (V^{-1})^*P^H_{H_0} J u \right)_H = 0,
\end{equation}
and
\begin{equation}
\label{f5_15}
\overline{
\left[ (V^{-1})^*P^H_{H_0} J u, (V^{-1})^*P^H_{H_0} J u \right]_J
}
\not= [u,u]_J.
\end{equation}
However, we do not know, whether such a case can ever happen.
\end{rmr}

\textbf{Proof of Theorem~\ref{t5_2}. }
By Theorem~\ref{t5_1}, the given operator $V$ can be extended to a $J$-isometric operator, defined on the whole $H$, if and only if
there exists a $J$-isometric operator $W_1$ in $H$, with the domain $H_1$, such that relation~(\ref{f5_9}) holds.
Let $u$ be a non-zero element in $H_1$, and $v$ be a non-zero element in $H\ominus JR(V)$.
An arbitrary linear operator $W_1$ in $H$, with the domain $H_1$, such that relation~(\ref{f5_9}) holds,
has the following form:
\begin{equation}
\label{f5_16}
W_1 (\beta u) = \beta \lambda v + \beta J(V^{-1})^*P^H_{H_0} J u,\qquad \beta\in \mathbb{C},
\end{equation}
where $\lambda\in \mathbb{C}$. On the other hand, an operator $W_1$ of the form~(\ref{f5_16}) with an arbitrary complex parameter $\lambda$
is a linear operator in $H$, with the domain $H_1$, such that relation~(\ref{f5_9}) holds.

It is easy to check that an operator $W_1$ of the form~(\ref{f5_16}) is $J$-isometric if and only if relation~(\ref{f5_11}) holds.
The statement about the correspondence~(\ref{f5_12}) follows from the formulas~(\ref{f5_10}) and~(\ref{f5_16}).
Let us check the last statement of the theorem.
An arbitrary $J$-isometric operator $W\supseteq V$, $D(W)=H$, is invertible~\cite[p. 18]{cit_300_Z}.
Moreover, we have $W^{-1}\subseteq JW^*J$. Since $H$ is finite-dimensional, then $WH=H$.
Therefore $W^{-1} = JW^*J$.
$\Box$

\begin{example} ($J$-isometric operator which has no $J$-unitary extensions)
\label{e5_1}
Let $H=\mathbb{C}^2$,
$\vec e_0 = \left( \begin{array}{cc} 1\\
0 \end{array} \right)$,
$\vec e_1 = \left( \begin{array}{cc} 0\\
1 \end{array} \right)$, and
\begin{equation}
\label{f5_18}
J (\alpha \vec e_0 + \beta \vec e_1) = \overline{\alpha} \vec e_1 + \overline{\beta} \vec e_0,\qquad \alpha, \beta\in \mathbb{C}.
\end{equation}
Observe that $J$ is a conjugation in $H$.
Consider the following operator $V$:
$$ V \alpha \vec e_0 = 0,\qquad \alpha\in \mathbb{C}, $$
with the domain $D(V) = \mathop{\rm Lin}\nolimits \{ \vec e_0 \}$.
The operator $V$ is $J$-isometric. Since it is not invertible, it has no $J$-unitary extensions.
\end{example}

\begin{example} ($J$-isometric operator which has a unique $J$-unitary extension inside the original Hilbert space)
\label{e5_2}
Let $H$, $\vec e_0$, $\vec e_1$, and $J$ be the same as in Example~\ref{e5_1}.
Consider the following operator $V$:
$$ V \alpha \vec e_0 = \alpha \vec e_0,\qquad \alpha\in \mathbb{C}, $$
with the domain $D(V) = \mathop{\rm Lin}\nolimits \{ \vec e_0 \} =: H_0$.
The operator $V$ is $J$-isometric.
Observe that $H\ominus JR(V) = H_0$.
Set $u=\vec e_1$, $v = \vec e_0$. Then $[u,u]_J = [v,v]_J = 0$, and
$$ (V^{-1})^*P^H_{H_0} J u = \vec e_0. $$
Equation~(\ref{f5_11}) takes the following form: $2\lambda = 0$.
Therefore $W = E_H$ is the unique $J$-unitary extension inside $H$.

\end{example}

\begin{example} ($J$-isometric operator which has exactly two $J$-unitary extensions inside the original Hilbert space)
\label{e5_3}
Let $H$, $\vec e_0$, $\vec e_1$  be the same as in Example~\ref{e5_1}. Set
\begin{equation}
\label{f5_20}
J (\alpha \vec e_0 + \beta \vec e_1) = \overline{\alpha} \vec e_0 + \overline{\beta} \vec e_1,\qquad \alpha, \beta\in \mathbb{C}.
\end{equation}
Consider the following operator $V$:
$$ V \alpha \vec e_0 = \alpha \vec e_0,\qquad \alpha\in \mathbb{C}, $$
with the domain $D(V) = \mathop{\rm Lin}\nolimits \{ \vec e_0 \} =: H_0$.
The operator $V$ is $J$-isometric.
Observe that $H\ominus JR(V) = \mathop{\rm Lin}\nolimits \{ \vec e_1 \} =: H_1$.
Set $u=v=\vec e_1$. Then $[u,u]_J = [v,v]_J = 1$, and
$$ (V^{-1})^*P^H_{H_0} J u = 0. $$
Equation~(\ref{f5_11}) takes the following form: $\lambda^2 -1 = 0$.
Therefore we have two possible $J$-unitary extensions of $V$ inside $H$.

\end{example}

We can apply Theorem~\ref{t5_2} to obtain some sufficient results for the solvability of the moment problem~C.

\begin{thm}
\label{t5_3}
Let the moment problem C be given with a set $\{ s_{n,m} \}_{n,m=0}^d$, $d\in \mathbb{N}$.
Suppose that
the matrix $\Gamma := ( s_{n,m} )_{n,m=0}^d$ is complex symmetric,
relation~(\ref{f3_19}) holds,
and
$\det\Gamma \not= 0$.
Consider the corresponding $J$-isometric operator $A_0$ from~(\ref{f4_5})
(acting in $\mathbb{C}^\rho$, $\rho = d+1$, with $\Lambda$ in~(\ref{f4_1}), provided by the corollary from Takagi's factorization).
Let $u$ be a non-zero element in $H\ominus D(A_0)$, and $v$ be a non-zero element in $H\ominus JR(A_0)$; $H_0 := D(A_0)$.
Suppose that equation~(\ref{f5_11})  with respect to an unknown $\lambda\in \mathbb{C}$, where $V=A_0$, has a solution.
Then the moment problem has a solution.

\end{thm}
\textbf{Proof.}
Observe that $\det\Lambda\not=0$. The Gram matrix
$G = ((x_n,x_m)_H)_{n,m=0}^d$ is equal to $\Lambda\Lambda^*$.
Then $\det G\not=0$, and $\{ x_n \}_{n=0}^\infty$ are linearly independent.
Thus, we have $\dim (H\ominus D(A_0)) = 1$, and $A_0$ is invertible.
Applying Corollary~\ref{c4_1} and Theorem~\ref{t5_2} we complete the proof.
$\Box$



\begin{center}
{\large\bf Truncated moment problems for $J$-self-adjoint, $J$-skew-self-adjoint and $J$-unitary operators.}
\end{center}
\begin{center}
{\bf S.M. Zagorodnyuk}
\end{center}

In this paper we study truncated moment problems for $J$-self-adjoint, $J$-skew-self-adjoint and $J$-unitary operators.
Conditions of the solvability are given. Some canonical solutions of the moment problems are constructed. As a by-product, some extension results for
$J$-skew-symmetric and $J$-isometric operators are obtained.

}

\end{document}